\documentclass[12pt]{article}

        \usepackage{amsfonts}
        \usepackage{amssymb}
        
        \def\pf{\medbreak\noindent{\bf Proof:}\enspace}

        \def\iff{\Longleftrightarrow}

        \title{A Pair of Optimal Inequalities Related to the Error Function}

        \author {Mary Beth Ruskai \thanks{supported by National Science
        Foundation Grant DMS-97-06981} \\ Department of
        Mathematics \\ University of Massachusetts  Lowell \\ Lowell,
        MA  01854 USA \\ {\normalsize bruskai@cs.uml.edu}
        \and Elisabeth Werner \\Department of Mathematics \\
    Case Western Reserve University \\
    Cleveland, OH 44106 USA \\
    {\normalsize emw2@po.cwru.edu}}

\begin{document}

        \maketitle

        \bigskip
        The function
        \begin{eqnarray}
        V(x)  & \equiv & \sqrt{\pi} e^{x^2} [1 - \hbox{erf}(x)]  \\
          & = &  \int_0^\infty  \frac{ e^{-u} }{\sqrt{x^2 + u}} du
                  = 2 e^{x^2}\int_x^\infty  e^{-t^2} dt \nonumber
        \end{eqnarray}
        arises in many contexts, from probability to mathematical physics,
        and satisifies the differential equation
        \begin{eqnarray}\label{Vdiffeq}
        V'(x) = 2x V(x) -2 ~~\hbox{with}~~  V(0) = \sqrt{\pi}.
        \end{eqnarray}
        In this note, we restrict attention to $x \geq 0$, in which case
        $ 0 < V(x) < 1/x$ and $V(x)$ is decreasing.  For $x \geq 0$ we show
that
        \begin{eqnarray}\label{opt.ineq}
          g_{\pi}(x) \leq  V(x) < g_4(x)
        \end{eqnarray}
        where
        \begin{eqnarray}\label{gkdef}
         g_k(x) = \frac{k}{(k-1)x + \sqrt{x^2+k}}.
        \end{eqnarray}
        We also show that
         these inequalities are optimal for functions of the form (\ref{gkdef})
        with equality only at $g_{\pi}(0) = V(0) = \sqrt{\pi}.$
        The bounds in (\ref{opt.ineq}) are considerably sharper than
        the classical inequalities of Komatsu
        \cite{IM}.  The upper bound, which implies the convexity of $1/V(x)$,
         was established independently by
        Wirth \cite{W} for $x \geq 0$ and by Szarek and Werner \cite{SzW} for
        $x > -\frac{1}{\sqrt 2}$.
        The weaker lower bound $g_3(x) < V(x)$ was used in \cite{BR} to show
        that the function $ \left[ 1/V(x) - x \right]^2/V(x)$ is decreasing
        for $x \geq 0$.  It is easy to see that the family of functions
        $g_k(x)$ is increasing in $k$ and that $ 0 < g_k(x) < 1/x.$

        In order to prove that the upper bound is optimal, we first observe
        that
        \begin{eqnarray}\label{ineq:deriv}
        \lefteqn{ g_k'(x) > 2[x g_k(x)  - 1] }  \nonumber \\
          & \iff & (k-2) x^2 + k (k-3) < (k-2) x \sqrt{x^2+k} \nonumber
        \\ & \iff & x^2 (k-2)(k-4) + k(k-3)^2 < 0,
        \end{eqnarray}
        when $k > 3.$  We now restrict attention to $3 \leq k \leq 4$ and
        let $h_k(x) = g_k(x) - V(x).$  For $k=4$, the expression
        (\ref{ineq:deriv}) implies  $g_4'(x) < 2[xg_k(x)  - 1]$
        so that $h_4^{\prime}(x) < 2x h_4(x)$ for all $ x \geq 0$; whereas for
        $k < 4$ this holds only for
        $x < a_k = \sqrt{\frac{k(k-3)^2}{(k-2)(4-k)}}.$
        Since both $V(x)$ and $g_k(x)$ are positive and bounded above
        by $1/x$, their difference also satisfies
        $|h_k(x)| < 1/x \rightarrow 0.$

        For $k=4$, if $h_4(x) \leq 0$ for some $x > 0 $, then  $h_4'(x)< 2x
h_4(x)$ is negative
        and thus $h_4$ is negative and strictly decreasing from a certain
$x$ on,
        which contradicts
        $\lim_{x \rightarrow \infty} h_4(x) = 0.$  Thus $h_4(x) > 0$ so
        that $g_4(x) > V(x)$, for all $x$.
        Now suppose that
        for some $k < 4$, $g_k$ is an upper bound, i.e. $h_k(x) \geq 0$
        for all $x \geq 0$.
        In particular, $h_k(x) \geq 0$ for all $x > a_k$.
        For $k < 4$, we find however
        that  $h_k^{\prime}(x) > 2 x h_k(x)$ holds for $x > a_k$.
 Thus we get $h_k(x) \geq 0$ and strictly increasing for all $x > a_k$
which
contradicts $\lim_{x \rightarrow \infty} h_k(x) = 0.$
        Thus the upper bound  can {\em not} hold when
         $x > a_k$ and $k < 4$. The lower bound also fails
         for $k > \pi$ since then
         $h_k(0) = g_k(0) - V(0) = \sqrt{k} - \sqrt{\pi} > 0.$

    To establish the improved lower bound $g_{\pi} \leq V(x)$
    we note that the argument above implies that
    $h_k(x)$ is negative for $x > a_k$ and $3 < k \leq \pi.$
    However for $k < \pi$ we have $h_k(0) < 0$ so that
    $h_k(x)$ is also negative for very small $x$.  If $h_k(x)$
     is ever non-negative, we can let $b$ denote the first
     place $h_k(x)$ touches or crosses the x-axis, i.e., $h_k(b) = 0$
     and $h_k(x) < 0$ for $x < b$.   Then $h_k$ must be increasing
     on some interval of the form $(x_0,b)$.  However, by the
     remarks above, $h_k(b) = 0$ implies $b \leq a_k$ so that
     $h_k^{\prime}(x) < 2 x h_k(x) < 0$ on $(x_0,b).$  Since this
     contradicts $h_k$ increasing on $(x_0,b)$, we must have
     $h_k(x) < 0$ for all $x \geq 0$ if $k < \pi.$

    Thus we have proved the lower bound
    $g_k(x) < V(x)$ on $[0,\infty)$  for $k < \pi.$
    Since $g_k$ is  continuous and increasing in $k$, it follows that
    $g_{\pi}(x) \leq V(x).$  To show that this inequality is strict
    except at $x=0$, note that the right derivative of $h_k$ at $0$
    satisfies $h_k'(0) = 2 - k$ so that  $h_{\pi}'(0) < 0$ and
    $h_{\pi}(x)$ is negative at least on some small interval $(0,x_1).$
    Then we can repeat the argument above to show that $h_{\pi}(x) < 0$
    if $x > 0.$
         \bigskip

        \noindent{\bf Acknowledgment} We are grateful for the hospitality
         of the School of Mathematics at Georgia Tech where this work was
         done.


\begin{thebibliography}{~~}


        \bibitem{IM} K. Ito and H.P. McKean {\em Diffusion Processes and
        Their Sample Paths} (Springer-Verlag, 1965).

        \bibitem{BR} R. Brummelhuis and M. B. Ruskai,
        ``A One-Dimensional Model for Many-Electron Atoms
        in Extremely Strong Magnetic Fields:  Maximum Negative Ionization''
            preprint.

        \bibitem{SzW} S. J. Szarek and E. Werner, ``Confidence Regions for
Means
        of Multivariate Normal Distributions and a Non-symmetric Correlation
        Inequality for Gaussian Measure'' MSRI preprint No. 1997-009,
        to appear in {\em J. Multivariate Analysis}.

        \bibitem{W} M. Wirth, ``On consid\`ere la fonction de ${\mathbb R}$
dans
        ${\mathbb R}$ d\'efine par $f(x) = e^{-x^2/2}$; d\'emontrer que la
        fonction $g$ de ${\mathbb R}$ dans ${\mathbb R}$ d\'efine par
        $g(x) = f(x)/\int_x^\infty f(t) dt$ est convexe''
        {\em Revue de Math\'ematiques Sp\'eciales} {\bf 104}, 187-88 (1993).


        \end{thebibliography}
         \end{document}